\documentclass[a4paper,12pt,reqno]{amsart}
\usepackage{amssymb}
\usepackage{amsmath}
\usepackage{amsthm}
\usepackage{amsfonts}
\usepackage{mathrsfs}

\newtheorem{theorem}{Theorem}[section]

\newtheorem{cor}[theorem]{Corollary}

\theoremstyle{definition}
\newtheorem{definition}[theorem]{Definition}

\theoremstyle{remark}

\numberwithin{equation}{section}

\allowdisplaybreaks[1]

\begin{document}

\title[Sampling theorem]{Sampling theorem and reconstruction formula for the space of translates on the Heisenberg group}

\author{S. Arati}
\author{R. Radha*}
\address{Department of Mathematics, Indian Institute of Technology Madras, Chennai 600 036, India}
\email{aratishashi@gmail.com ; radharam@iitm.ac.in}
\thanks{* Corresponding author}



\begin{abstract}
The paper deals with the necessary and sufficient conditions for obtaining reconstruction formulae and sampling theorems for every function belonging to the principal shift invariant subspace of $L^2(\mathbb{H}^n)$, both in the time domain and a transform domain, where $\mathbb{H}^n$ denotes the Heisenberg group.
\end{abstract}

\maketitle
\section{Introduction}
The classical Shannon sampling theorem \cite{Shan} states that if $f\in L^{2}(\mathbb{R})$ whose Fourier transform has support in $[-\frac{1}{2},\frac{1}{2}]$, then $f$ can be recovered from its uniform samples at integers by the formula
\[f(x)=\sum\limits_{k\in\mathbb{Z}}f(k)\dfrac{\sin\pi (x-k)}{\pi (x-k)}.\]
The reconstruction of such band limited functions $f$ from nonuniform samples dates back to the works of Paley and Wiener \cite{Paley} and Duffin and Eachus \cite{Duff_Each} in Mathematics literature. Later, Kadec \cite{Kadec} proved that in order to obtain a stable set of sampling for such functions from the nonuniform samples the perturbation should be at the most $\tfrac{1}{4}$ from the integers. The more general sampling condition was given in terms of Beurling density in \cite{Land}. We refer to the work of Butzer and Stens \cite{Butz} for the classical historical review of sampling theory. 
\par
In \cite{Walt}, Walter extended the Shannon sampling theorem to wavelet subspaces. In 1999, Zhou and Sun \cite{Zhou_Sun} provided a necessary and sufficient condition under which every function in a closed subspace $V_0$ of $L^2(\mathbb{R})$ has a sampling expansion. In other words, they proved necessary and sufficient conditions to obtain the sampling formula $f(x)=\sum_{k}f(k)g(x-k)$, $\forall\,f\in V_0$ which holds for some $\{g(\cdot-n):n\in\mathbb{Z}\}$ in $V_0$ with the convergence being both in $L^2(\mathbb{R})$ and pointwise on $\mathbb{R}$. 
\par
The problem of sampling and reconstruction has been studied in a shift invariant space with a single generator by several authors. (See for example \cite{Liu, Ald_Groch, SunW_Zhou, Garci2, Anto_Rad}.) For the sampling problem in a shift invariant space with multiple generators, we refer to \cite{Ald_Sun_Tang, Ald_Krish, Sun1, Xian_Li, Garci, Acost_Ald_Krish, Nash_Sun, Sun2, Nash_Sun_Xian, Zhang, Fuhr_Xian, Rad_Sarv_Siva}. 
\par
In \cite{RaSiva}, a Shannon type sampling theorem was proved for the Heisenberg group. Recently, Radha and Saswata obtained a sampling theorem on a subspace of a twisted shift-invariant space in \cite{RaS4}. In fact, they gave a necessary and sufficient condition for obtaining a reconstruction formula for functions belonging to a subspace $V^{0,t}(\varphi)=\overline{span}\{T^t_{(k,0)}\varphi:k\in\mathbb{Z}^n\}$ of the principal twisted shift invariant space $V^{t}(\varphi)=\overline{span}\{T^t_{(k,l)}\varphi:(k,l)\in\mathbb{Z}^{2n}\}$ in $L^2(\mathbb{R}^{2n})$ from their samples $\{f(k,j):k\in\mathbb{Z}^n\}$ for each fixed $j\in\mathbb{Z}^n$.
\par
The aim of our paper is to look for necessary and sufficient conditions for obtaining a sampling theorem and a reconstruction formula for every function belonging to the principal shift invariant space of $L^2(\mathbb{H}^n)$. We also consider a transform $F_{\lambda}$ of functions $f$ in $L^2(\mathbb{H}^n)$ which is an operator valued function on $\mathbb{R}^{2n}$ and derive reconstruction formula in the transform domain of $f$. We provide an inversion theorem which recovers $f$ from $F_{\lambda}$, thereby paving way for another reconstruction formula for $f$ in the time domain. We may also view the reconstruction formulae as matrix equations which aids in getting corresponding sampling theorems.

\section{Necessary background}
Let $\mathcal{H}$ be a separable Hilbert space.
\begin{definition}
A sequence $\{f_{k}:k\in\mathbb{Z}\}$ in $\mathcal{H}$ is a frame for $\mathcal{H}$ if there exist constants $A,B>0$ such that \[A\|f\|^{2}\leq\sum_{k\in\mathbb{Z}}|\langle f,f_k\rangle|^{2}\leq B\|f\|^{2},\quad \forall\, f\in\mathcal{H}.\]The numbers $A$ and $B$ are called frame bounds. If the right hand side inequality holds, then $\{f_{k}:k\in\mathbb{Z}\}$ is said to be a Bessel sequence with bound $B$. A sequence $\{f_k:k\in\mathbb{Z}\}$ in $\mathcal{H}$ is said to be a frame sequence if it is a frame for $\overline{\text{span}}\{f_k:k\in\mathbb{Z}\}$.
\end{definition}

\begin{definition}
 A Riesz basis for $\mathcal{H}$ is a family of the form $\{Ue_{k}:k\in\mathbb{Z}\}$, where $\{e_{k}:k\in\mathbb{Z}\}$ is an orthonormal basis for $\mathcal{H}$ and $U:\mathcal{H}\rightarrow\mathcal{H}$ is a bounded invertible operator. Alternatively, a sequence $\{f_{k}:k\in\mathbb{Z}\}$ is a Riesz basis for $\mathcal{H}$ if $\{f_{k}:k\in\mathbb{Z}\}$ is complete in $\mathcal{H}$, and there exist constants $A,B>0$ such that for every finite scalar sequence $\{c_k\}$, one has 
\[A\sum_{k}|c_k|^2\leq \|\sum_{k}c_k f_k\|^2\leq B\sum_{k}|c_k|^2.\]
A sequence $\{f_k:k\in\mathbb{Z}\}$ in $\mathcal{H}$ is a Riesz sequence if it is a Riesz basis  for $\overline{\text{span}}\{f_k:k\in\mathbb{Z}\}$.    
\end{definition}

\indent  Let $\{f_{k}:k\in\mathbb{Z}\}$ be a frame for $\mathcal{H}$. The operator \[S:\mathcal{H}\rightarrow\mathcal{H},\quad Sf=\sum_{k\in\mathbb{Z}}\langle f,f_{k}\rangle f_{k}\] is called the frame operator. It is bounded, invertible, self-adjoint and positive.

\begin{definition}
Suppose $\{f_{k}:k\in\mathbb{Z}\}$ is a frame for $\mathcal{H}$. The canonical dual frame of $\{f_{k}:k\in\mathbb{Z}\}$ is the frame $\{S^{-1}f_{k}:k\in\mathbb{Z}\}$ for $\mathcal{H}$, where $S$ is the frame operator.
\end{definition}

Any element $f\in\mathcal{H}$ can be expressed in terms of $S^{-1}f_{k}$ as follows.
\begin{enumerate}
\item[(i)] $f=\sum_{k\in\mathbb{Z}}\langle f, S^{-1}f_{k}\rangle f_k$,
\item[(ii)] $f=\sum_{k\in\mathbb{Z}}\langle f, f_{k}\rangle S^{-1}f_k$.
\end{enumerate}

\begin{definition}
Let $\{f_{k}:k\in\mathbb{Z}\}$ be a frame for $\mathcal{H}$. A frame $\{g_{k}:k\in\mathbb{Z}\}$ in $\mathcal{H}$ such that 
\[f=\sum_{k\in\mathbb{Z}}\langle f, g_k\rangle f_k,\quad\forall f\in\mathcal{H}\]
is called a dual frame of $\{f_{k}:k\in\mathbb{Z}\}$.
\end{definition}

\indent
Let $\{f_{k}:k\in\mathbb{Z}\}$ be a frame for $\mathcal{H}$ which is not a Riesz basis. Then there exist dual frames other than the canonical dual frame. So, every element
in $\mathcal{H}$ has a representation in terms of the frame
elements, namely, 
\[f=\sum_{k\in\mathbb{Z}}c_kf_k,\quad f\in\mathcal{H},\]
where $\{c_k\}\in l^2(\mathbb{Z})$. The above representation is called a frame expansion of $f$.

\par
For more on frames and bases, we refer to the books by Christensen\cite{Chris} and Heil\cite{Heil}.\\
\par
The Heisenberg group $\mathbb{H}^n$ is a nilpotent Lie group whose underlying manifold is $\mathbb{R}^n\times\mathbb{R}^n\times\mathbb{R}$ which satisfies the group law
\[(x,y,t)(u,v,s)=\left(x+u,y+v,t+s+\frac{1}{2}(u\cdot y-v\cdot x)\right).\]
It is a nonabelian noncompact locally compact group. The Haar measure on $\mathbb{H}^n$ is the Lebesgue measure $dxdydt$. From the well known Stone-von Neumann theorem it follows that every infinite dimensional irreducible unitary representation on the Heisenberg group is unitarily equivalent to the representation $\pi_\lambda,\,\lambda\in\mathbb{R}^\ast$, where $\pi_\lambda$ is defined by
\[\pi_\lambda(x,y,t)\varphi(\xi)=e^{2\pi i\lambda t}e^{2\pi i\lambda(x\cdot\xi+\frac{1}{2}x\cdot y)}\varphi(\xi+y),\quad\varphi\in L^{2}(\mathbb{R}^n).\]
For $f\in L^{1}(\mathbb{H}^n)$, the group Fourier transform $\hat{f}$ is defined as follows. For $\lambda\in\mathbb{R}^\ast,\,\hat{f}(\lambda)$ is given by 
\[\hat{f}(\lambda)=\int\limits_{\mathbb{C}^n\times\mathbb{R}}f(z,t)\pi_\lambda(z,t)dzdt.\] More explicitly, $\hat{f}(\lambda)$ is the bounded operator acting on $L^{2}(\mathbb{R}^n)$ (i.e., $\hat{f}(\lambda)\in\mathscr{B}(L^2(\mathbb{R}^n))$) given by $\hat{f}(\lambda)\varphi=\int\limits_{\mathbb{C}^n\times\mathbb{R}}f(z,t)\pi_\lambda(z,t)\varphi dzdt,\,\varphi\in L^{2}(\mathbb{R}^n)$, where the integral is a Bochner integral taking values in the Hilbert space $L^{2}(\mathbb{R}^n)$. Further, $\|\hat{f}(\lambda)\|_{\mathscr{B}}\leq\|f\|_{L^{1}(\mathbb{H}^n)}$.
The inverse Fourier transform of $f\in L^{1}(\mathbb{H}^n)$ in the $t$ variable, denoted by $f^\lambda$, is defined as  
\begin{equation*}\label{E:f_lambda}
f^\lambda(z)=\int\limits_{\mathbb{R}}f(z,t)e^{2\pi i\lambda t}dt.
\end{equation*}
It can be seen that $f^\lambda\in L^{1}(\mathbb{C}^n)$. For $f\in L^{1}(\mathbb{C}^n)$, the operator $W_\lambda(f)$ on $L^{2}(\mathbb{R}^n)$ is defined as \[W_\lambda(f)=\int\limits_{\mathbb{C}^n}f(z)\pi_\lambda(z,0)dz.\] Clearly, there is a relation between group Fourier transform and $W_\lambda$ given by
\begin{equation*}\label{E:Reln bet FT and W}
\hat{f}(\lambda)=W_\lambda(f^\lambda).
\end{equation*}
Moreover, $W_\lambda(f)$ is an integral operator on $L^{2}(\mathbb{R}^n)$ with kernel $K^\lambda_f$ given by
\begin{equation*}
K^\lambda_f(\xi,\eta)=\int\limits_{\mathbb{R}^n}f(x,\eta-\xi)e^{\pi i\lambda x\cdot(\xi+\eta)}dx.
\end{equation*} 
In particular when $\lambda=1,\,W_\lambda(f)$ is denoted by $W(f)$ which is called the Weyl transform of $f$ and the associated kernel is denoted by $K_f$.
\par
If $f$ and $g$ are in $L^1(\mathbb{H}^n)$, then their convolution is defined by
\[f\ast g(z,t)=\int\limits_{\mathbb{C}^n\times\mathbb{R}}f((z,t)(-w,-s))g(w,s)dwds.\]
Under this convolution, $L^{1}(\mathbb{H}^n)$ turns out to be a noncommutative Banach algebra. The group Fourier transform takes convolution into products like in the classical case, i.e., \[(f\ast g\hat{)}(\lambda)=\hat{f}(\lambda)\hat{g}(\lambda).\]

\par
Analogous to the case of the Euclidean Fourier transform, the definitions of $W_\lambda$ and the group Fourier transform $\hat{f}$ can be extended to functions in $L^{2}(\mathbb{C}^n)$ and $L^{2}(\mathbb{H}^n)$ respectively through the density argument. In fact, for $f\in L^{2}(\mathbb{C}^n),\,W_\lambda(f)$ is a Hilbert-Schmidt operator on $L^{2}(\mathbb{R}^n)$ which satisfies
\[\|W_\lambda(f)\|_{\mathcal{B}_2}=\|K^\lambda_f\|_{L^{2}(\mathbb{C}^n)}=\tfrac{1}{|\lambda|^{n/2}}\|f\|_{L^{2}(\mathbb{C}^n)},\]
where $\mathcal{B}_2=\mathcal{B}_2(L^{2}(\mathbb{R}^n))$ denotes the class of Hilbert-Schmidt operators on $L^{2}(\mathbb{R}^n)$. In other words, for $f,g\in L^{2}(\mathbb{C}^n)$,
\begin{equation*}\label{E:Planch for Weyl T}
\langle W_\lambda(f),W_\lambda(g)\rangle_{\mathcal{B}_2}=\langle K^\lambda_f,K^\lambda_g\rangle_{L^{2}(\mathbb{C}^n)}=\tfrac{1}{|\lambda|^n}\langle f,g\rangle_{L^{2}(\mathbb{C}^n)}.
\end{equation*}
Furthermore, the group Fourier transform satisfies the Plancherel formula 
\[\|\hat{f}\|_{L^{2}(\mathbb{R}^\ast,\mathcal{B}_2;d\mu)}=\|f\|_{L^{2}(\mathbb{H}^n)},\]
where $L^{2}(\mathbb{R}^\ast,\mathcal{B}_2;d\mu)$ stands for the space of functions on $\mathbb{R}^\ast$ taking values in $\mathcal{B}_2$ and square integrable with respect to the measure $d\mu(\lambda)=|\lambda|^n d\lambda.$ 
\par
For further study on the Heisenberg group, we refer to \cite{Foll} and \cite{Thanga}. 
\par
For $f\in L^2(\mathbb{H}^n)$ and $\lambda\in\mathbb{R}^\ast$, an operator valued function, $F_{\lambda}$,  on $\mathbb{R}^{2n}$ is defined as
\begin{equation}\label{E: F lambda}
F_{\lambda}(\xi_1,\xi_2)=\pi_{\lambda}(\xi_1,\xi_2,0)\hat{f}(\lambda)\pi_{\lambda}(-\xi_1,-\xi_2,0),\quad (\xi_1,\xi_2)\in\mathbb{R}^{2n}.
\end{equation}
The above operator was used by Thangavelu \cite{Thanga_Paley} in studying Paley-Wiener theorem for the Heisenberg group.\\
\par

Let $\mathscr{L}$ denote a lattice in $\mathbb{H}^{n}$. In other words, $\mathscr{L}$ is a discrete subgroup of the Heisenberg group $\mathbb{H}^{n}$ such that $\mathbb{H}^{n}/\mathscr{L}$ is compact. For $\varphi\in L^{2}(\mathbb{H}^{n})$, the shift invariant space, $V(\varphi)$, is defined to be $\overline{span}\{{L_{l}\varphi:l\in\mathscr{L}}\}$, where
$L_{l}\varphi(X)=\varphi(l^{-1}\cdot X),~X\in\mathbb{H}^{n}$. However, from the computational point of view, one can work with the standard lattice $\{(2k,l,m):k,l\in\mathbb{Z}^{n},m\in\mathbb{Z}\}$ in place of $\mathscr{L}$. Explicitly, the action of the left translation $L_{(2k,l,m)}$ on $L^2(\mathbb{H}^n)$ for $(2k,l,m)\in\mathbb{Z}^n\times\mathbb{Z}^n\times\mathbb{Z}$ is given by
\begin{align*}
L_{(2k,l,m)}\varphi(x,y,t)&=\varphi((2k,l,m)^{-1}(x,y,t))\\
&=\varphi\left(x-2k,y-l,t-m+\tfrac{1}{2}(y\cdot 2k-x\cdot l)\right), 
\end{align*}
for $\varphi\in L^2(\mathbb{H}^n)$ and the principal shift invariant space generated by $\varphi$ is given by
\[V(\varphi)=\overline{span}\{L_{(2k,l,m)}\varphi:(k,l,m)\in\mathbb{Z}^n\times\mathbb{Z}^n\times\mathbb{Z}\}.\]
The left translates satisfy the following equations.
\begin{align*}
&\widehat{L_{(2k,l,m)}\varphi}(\lambda)=e^{2\pi im\lambda}\widehat{L_{(2k,l,0)}\varphi}(\lambda)\\
&\|\widehat{L_{(2k,l,0)}\varphi}(\lambda)\|_{\mathcal{B}_2}=\|\widehat{\varphi}(\lambda)\|_{\mathcal{B}_2},\label{E: norm of FT of trans}  
\end{align*}
for $\lambda\in\mathbb{R}^{\ast}$.
For a study of frames and Riesz bases in connection with shift invariant spaces on $\mathbb{H}^n$, we refer to \cite{Azit, RaS2}. In \cite{RaS2}, the system of left translates on $\mathbb{H}^n$ has been characterized to be a frame sequence and a Riesz sequence in terms of the weight function $G^{\varphi}_{k,l}$ defined below.
\begin{definition}
For $\varphi\in L^2(\mathbb{H}^n)$ and $k,l\in\mathbb{Z}^{n}$, the function $G^{\varphi}_{k,l}$ is defined by
\begin{equation*}\label{D: G kl}
G_{k,l}^{\varphi}(\lambda)=\sum\limits_{r\in\mathbb{Z}}\langle\widehat{\varphi}
(\lambda+r),\widehat{L_{(2k,l,0)}\varphi}(\lambda+r)\rangle_{\mathcal{B}_{2}}
|\lambda+r|^{n},\;\lambda\in (0,1].
\end{equation*}
\end{definition}
It can also be written in terms of the kernel of $W_{\lambda}$ as 
\begin{align*}
G^\varphi_{k,l}(\lambda)=\sum\limits_{r\in\mathbb{Z}}\int\limits_{\mathbb{R}^n}\int\limits_{\mathbb{R}^n}&
K^{\lambda+r}_{\varphi^{\lambda+r}}(\xi,\eta)\overline{K^{\lambda+r}_{\varphi^{\lambda+r}}(\xi+l,\eta)}e^{-2\pi i(\lambda+r)k\cdot(2\xi+l)}d\xi d\eta|\lambda+r|^n.
\end{align*}
We refer to \cite{RaS2} in this connection.

\begin{definition}[\cite{RaS2}]
A function $\varphi\in L^{2}(\mathbb{H}^{n})$ is said to satisfy Condition C if $G_{k,l}^{\varphi}(\lambda)=0$ a.e. $\lambda\in (0,1]$,~for all $(k,l)\in\mathbb{Z}^{2n}\backslash\{(0,0)\}$.
\end{definition}

The characterizations for the system of left translates on $\mathbb{H}^n$ to be a Bessel sequence and a frame sequence are as follows.

\begin{theorem}[\cite{RaS2}]\label{T:Bessel char}
Let $\varphi\in L^2(\mathbb{H}^n)$. Suppose $\{L_{(2k,l,m)}\varphi:(k,l,m)\in\mathbb{Z}^n\times\mathbb{Z}^n\times\mathbb{Z}\}$ is a Bessel sequence in $L^{2}(\mathbb{H}^{n})$ with bound $B$. Then $G^{\varphi}_{0,0}(\lambda)\leq B\quad \text{a.e. }\lambda\in (0,1]$. Conversely, suppose there exists $B>0$ such that $G^{\varphi}_{0,0}(\lambda)\leq B\quad \text{a.e. }\lambda\in (0,1]$ and $\varphi$ satisfies Condition $C$. Then $\{L_{(2k,l,m)}\varphi:(k,l,m)\in\mathbb{Z}^n\times\mathbb{Z}^n\times\mathbb{Z}\}$ is a Bessel sequence in $L^{2}(\mathbb{H}^{n})$ with bound $B$.
\end{theorem}

\begin{theorem}[\cite{RaS2}]\label{T: frame char}
Let $\varphi\in L^2(\mathbb{H}^n)$ satisfy Condition C. Then the collection $\{L_{(2k,l,m)}\varphi:(k,l,m)\in\mathbb{Z}^n\times\mathbb{Z}^n\times\mathbb{Z}\}$ is a frame for $V(\varphi)$ with frame bounds $A,B>0$ if and only if
\[A\leq G^{\varphi}_{0,0}(\lambda)\leq B\quad \text{a.e. }\lambda\in\Omega_{\varphi},\]
where $\Omega_{\varphi}=\{\eta\in(0,1]: G^{\varphi}_{0,0}(\eta)>0\}.$
\end{theorem}

\section{The main results}
We shall now state some necessary conditions when a sampling formula holds for the principal shift invariant space $V(\varphi)$, $\varphi\in L^2(\mathbb{H}^n)$.

\begin{theorem}\label{T:Nec}
Let $V_0$ be a closed subspace of $L^2(\mathbb{H}^n)$. Let $\varphi\in L^2(\mathbb{H}^n)$ satisfy Condition C and the collection $\{L_{(2k,l,m)}\varphi:(k,l,m)\in\mathbb{Z}^n\times\mathbb{Z}^n\times\mathbb{Z}\}$ be a frame for $V_0$. Suppose, for each $(k,l)\in\mathbb{Z}^{2n}$, $\sum_{m\in\mathbb{Z}}c_m L_{(2k,l,m)}\varphi$ converges pointwise to a continuous function, for any $\{c_m\}\in l^2(\mathbb{Z})$. Further, suppose there exists a function $\psi\in L^2(\mathbb{H}^n)$ satisfying Condition C such that $\{L_{(2k,l,m)}\psi:(k,l,m)\in\mathbb{Z}^n\times\mathbb{Z}^n\times\mathbb{Z}\}$ is a frame for $V_0$ and 
\[f=\sum\limits_{(k,l,m)\in\mathbb{Z}^{2n+1}}f(2k,l,m)L_{(2k,l,m)}\psi,\quad\forall\, f\in V_0,\]
with the convergence in $L^2(\mathbb{H}^n)$. Then $\varphi\in C(\mathbb{H}^n)$, 
\[\sup_{(x,y,t)\in\mathbb{T}^{2n+1}}\sum_{m\in\mathbb{Z}}|L_{(2k,l,m)}\varphi(x,y,t)|^2<\infty,\]
for each $(k,l)\in\mathbb{Z}^{2n}$ and 
\[A\chi_{\Omega_{\varphi}}(\lambda)\leq \|P(\lambda)\|^2_{l^2(\mathbb{Z}^{2n})}\leq B\chi_{\Omega_{\varphi}}(\lambda)\quad \text{a.e. }\lambda\in\mathbb{R}\]
for some constants $A,B>0$, where
\begin{equation*}\label{E:P kl defn}
P(\lambda)=\{P_{k,l}(\lambda)\}_{(k,l)\in\mathbb{Z}^{2n}}\text{ with }P_{k,l}(\lambda)=\sum_{m\in\mathbb{Z}}\varphi(2k,l,m)e^{2\pi im\lambda},
\end{equation*}
$\Omega_{\varphi}=\{\lambda\in\mathbb{R}:G^{\varphi}_{0,0}(\lambda)>0\}$ and $G^{\varphi}_{0,0}$ is as in Definition \ref{D: G kl}.
\end{theorem}

The following theorem gives sufficient conditions for obtaining a reconstruction formula for the principal shift invariant subspace of $L^2(\mathbb{H}^n)$.

\begin{theorem}\label{T:suff condn in time domain}
Let $V_0$ be a closed subspace of $L^2(\mathbb{H}^n)$. Let $\varphi\in L^2(\mathbb{H}^n)$ satisfy Condition C and the collection $\{L_{(2k,l,m)}\varphi:(k,l,m)\in\mathbb{Z}^n\times\mathbb{Z}^n\times\mathbb{Z}\}$ be a frame for $V_0$. For $f\in V_0$, let
\[f=\sum_{(k,l,m)\in\mathbb{Z}^{2n+1}}c_{2k,l,m}L_{(2k,l,m)}\varphi
=\sum_{(k,l)\in\mathbb{Z}^{2n}}g_{2k,l},\]
where 
$g_{2k,l}=\sum_{m\in\mathbb{Z}}c_{2k,l,m}L_{(2k,l,m)}\varphi$
in $L^2(\mathbb{H}^n)$. Suppose $\varphi\in C(\mathbb{H}^n)$ and for each $(k,l)\in\mathbb{Z}^{2n}$, 
\[\sup_{(x,y,t)\in\mathbb{H}^{n}}\sum_{m\in\mathbb{Z}}|L_{(2k,l,m)}\varphi(x,y,t)|^2<\infty.\] Further, suppose there exist constants $A,B>0$ such that
\[A\chi_{\Omega_{\varphi}}(\lambda)\leq \|P(\lambda)\|^2_{l^2(\mathbb{Z}^{2n})}\leq B\chi_{\Omega_{\varphi}}(\lambda)\quad \text{a.e. }\lambda\in\mathbb{R},\]
where $P(\lambda)$ and $\Omega_{\varphi}$ are as in Theorem \ref{T:Nec}. Then, for each $(k,l)\in\mathbb{Z}^{2n}$, $\sum_{m\in\mathbb{Z}}c_m L_{(2k,l,m)}\varphi$ converges pointwise to a continuous function, for any $\{c_m\}\in l^2(\mathbb{Z})$ and there exists a function $\psi\in L^2(\mathbb{H}^n)$ satisfying Condition C such that $\{L_{(2k,l,m)}\psi:(k,l,m)\in\mathbb{Z}^n\times\mathbb{Z}^n\times\mathbb{Z}\}$ is a frame for $V_0$ with the following reconstruction formula. For every $f\in V_0$,  
\[f=\sum\limits_{(k,l,m)\in\mathbb{Z}^{2n+1}}\alpha_{2k,l,m}L_{(2k,l,m)}\psi,\]
with the convergence in $L^2(\mathbb{H}^n)$, where $\alpha_{2k,l,m}=g_{2k,l}\ast_d\tilde{\varphi}(2k,l,m)$, $\tilde{\varphi}(X)=\overline{\varphi(X^{-1})}$, $X\in\mathbb{H}^n$ and $g_{2k,l}\ast_d\tilde{\varphi}(2k',l',m')$ denotes $$\sum_{(k'',l'',m'')\in\mathbb{Z}^{2n+1}}g_{2k,l}((2k',l',m')(2k'',l'',m'')^{-1})\tilde{\varphi}(2k'',l'',m''),$$ 
for $(2k',l',m')\in\mathbb{Z}^{2n+1}.$
\end{theorem}

Next, we provide a sufficient condition for obtaining an exact sampling formula for a smaller class of functions, namely, $\{L_{(0,0,m)}\varphi:m\in\mathbb{Z}\}$. Let the closed linear span of this collection be denoted by $V^0(\varphi)$.

\begin{theorem}
Let $\varphi\in L^2(\mathbb{H}^n)$ be such that $\{L_{(0,0,m)}\varphi:m\in\mathbb{Z}\}$ is a frame for $V^0(\varphi)$. Suppose $\varphi\in C(\mathbb{H}^n)$ and 
\[\sup_{(x,y,t)\in\mathbb{H}^{n}}\sum_{m\in\mathbb{Z}}|L_{(0,0,m)}\varphi(x,y,t)|^2<\infty.\] Further, suppose there exist constants $A,B>0$ such that
\[A\chi_{\Omega_{\varphi}}(\lambda)\leq |P_0(\lambda)|\leq B\chi_{\Omega_{\varphi}}(\lambda)\quad \text{a.e. }\lambda\in\mathbb{R},\]
where $P_0(\lambda)=\sum_{m\in\mathbb{Z}}\varphi(0,0,m)e^{2\pi im\lambda}$ and $\Omega_{\varphi}$ is as in Theorem \ref{T:Nec}. Then, $\sum_{m\in\mathbb{Z}}c_m L_{(0,0,m)}\varphi$ converges pointwise to a continuous function, for any $\{c_m\}\in l^2(\mathbb{Z})$ and there exists a function $\psi\in L^2(\mathbb{H}^n)$ such that $\{L_{(0,0,m)}\psi:m\in\mathbb{Z}\}$ is a frame for $V^0(\varphi)$ with the following sampling formula. For every $f\in V^0(\varphi)$,  
\[f=\sum\limits_{m\in\mathbb{Z}}f(0,0,m)L_{(0,0,m)}\psi,\]
where the convergence is both in $L^2(\mathbb{H}^n)$ and uniform on $\mathbb{H}^n$.
\end{theorem}

Now, we shall obtain a sampling theorem for $f$, belonging to the principal shift invariant space of $L^2(\mathbb{H}^n)$, in the transform domain (Theorem \ref{T:F lambda sampling}) using the operator valued function $F_{\lambda}$ on $\mathbb{R}^{2n}$ given in (\ref{E: F lambda}). Using the inversion formula for $f$ from $F_{\lambda}$, we can also obtain a reconstruction formula for $f$ in the time domain (Corollary \ref{C: sampling for f using inversion}).  
\par
We note that for $(\xi_1,\xi_2)\in\mathbb{R}^{2n}$, $F_{\lambda}(\xi_1,\xi_2)$ can also be written as 
\begin{equation*}\label{E: alt exp for F lambda}
F_{\lambda}(\xi_1,\xi_2)=\int_{\mathbb{R}^{2n}}e^{2\pi i\lambda(x\cdot\xi_2-y\cdot\xi_1)}f^{\lambda}(x,y)\pi_{\lambda}(x,y,0)dxdy.
\end{equation*}
\begin{theorem}\label{T:F lambda sampling}
Let $\varphi\in L^2(\mathbb{H}^n)$ and $\{L_{(2k,l,m)}\varphi:(k,l,m)\in\mathbb{Z}^n\times\mathbb{Z}^n\times\mathbb{Z}\}$ be a frame sequence. For $f\in V(\varphi)$ and $\lambda\in\mathbb{R}^{\ast}$, the transform $F_{\lambda}$ of $f$, as defined in (\ref{E: F lambda}), satisfies
\[F_{\lambda}(\xi_1,\xi_2)=\sum_{k,l,m}c_{2k,l,m}e^{2\pi i\lambda(2k\cdot\xi_2-l\cdot\xi_1+m)}(\Phi^{M_{\lambda,\xi_1,\xi_2}})_{\lambda}(2k,l)\pi_{\lambda}(2k,l,0),\]
where $\{c_{2k,l,m}\}$ is the sequence of frame coefficients in a frame expansion of $f$ and $(\Phi^{M_{\lambda,\xi_1,\xi_2}})_{\lambda}$ is the corresponding transform of the modulated $\varphi$, namely $M_{(\lambda\xi_2,-\lambda\xi_1,0)}\varphi$, in $L^2(\mathbb{H}^n)$.  
\end{theorem}

The following is an inversion formula that gives $f\in L^2(\mathbb{H}^n)$ from $F_{\lambda}(\xi_1,\xi_2)$. 

\begin{theorem}\label{T:inversion formula}
For $f\in L^2(\mathbb{H}^n)$ and any $(\xi_1,\xi_2)\in\mathbb{R}^{2n}$, one has
\[\int_{\mathbb{R}}|\lambda|^n e^{2\pi i\lambda(v\cdot\xi_1-u\cdot\xi_2)}tr(\pi_\lambda(u,v,w)^{\ast}F_{\lambda}(\xi_1,\xi_2))d\lambda=f(u,v,w),\] where $(u,v,w)\in\mathbb{H}^n$. In particular,
\[\int_{\mathbb{R}}|\lambda|^n tr(\pi_\lambda(u,v,w)^{\ast}F_{\lambda}(0,0))d\lambda=f(u,v,w),\text{ for }(u,v,w)\in\mathbb{H}^n.\]
\end{theorem}

Using Theorems \ref{T:inversion formula} and \ref{T:F lambda sampling}, we obtain another reconstruction formula for $f$ in the time domain which is as follows.
\begin{cor}\label{C: sampling for f using inversion}
Let $\varphi\in L^2(\mathbb{H}^n)$ and $\{L_{(2k,l,m)}\varphi:(k,l,m)\in\mathbb{Z}^n\times\mathbb{Z}^n\times\mathbb{Z}\}$ be a frame sequence. For $f\in V(\varphi)$ and $(u,v,w)\in\mathbb{H}^n$,
\begin{align*}
f(u,v,w)=\int_{\mathbb{R}}&|\lambda|^n 
\sum_{(k,l,m)\in\mathbb{Z}^{2n+1}}c_{2k,l,m}e^{2\pi i\lambda m}
\\&\times tr\left(\pi_\lambda(u,v,w)^{\ast}\Phi_{\lambda}(2k,l)\pi_{\lambda}(2k,l,0) \right)d\lambda,
\end{align*}
where $\{c_{2k,l,m}\}$ is the sequence of frame coefficients in a frame expansion of $f$ and $\Phi_{\lambda}$ is the transform of $\varphi$, as defined in (\ref{E: F lambda}).
\end{cor}
    
\bibliographystyle{plain}
\bibliography{references}

\end{document}